\documentclass[12pt]{amsart}
\usepackage[all]{xy}
\usepackage{amsmath, amscd, graphicx, latexsym, hyperref, rlepsf, times}

\calclayout
\newtheorem{dummy}{anything}[section]
\newtheorem{theorem}[dummy]{Theorem}

\newtheorem{lemma}[dummy]{Lemma}

\theoremstyle{definition}

\newtheorem{remark}[dummy]{Remark}

\newcommand{\bR}{\mathbb R}
\newcommand{\tr}{\textrm}
\begin{document}
\title[Lutz Twist]
{Full Lutz twist along the binding of an open book}
\author{Burak Ozbagci and Mehmetcik Pamuk}
\address{Department of Mathematics \\ Ko\c{c} University \\ Istanbul, Turkey}
\email{bozbagci@ku.edu.tr}
\email{mpamuk{@}ku.edu.tr}
\subjclass[2000]{57R17, 57R65} \keywords{Lutz twist, contact
surgery, open book decomposition}
\begin{abstract}\noindent
Let $T$ denote a binding component of an open book $(\Sigma,
\phi)$ compatible with a closed contact $3$-manifold $(M, \xi)$.
We describe an explicit open book $(\Sigma', \phi')$ compatible
with $(M, \zeta)$, where $\zeta$ is the contact structure obtained
from $\xi$ by performing a full Lutz twist along $T$. Here,
$(\Sigma', \phi')$ is obtained from $(\Sigma, \phi)$ by a
\emph{local} modification near the binding.

\end{abstract}
\maketitle
\section{Introduction}

Let $T$ denote a binding component of an open book $(\Sigma,
\phi)$ compatible with a closed contact $3$-manifold $(M, \xi)$.
Then, by definition, $T$ is a transverse knot. By performing  a
full Lutz twist along $T$, we get a new contact structure $\zeta$
on $M$. Our intention in the present note is to give an explicit
open book $(\Sigma', \phi')$ compatible with $(M, \zeta)$.

Our construction can be outlined as follows. First we use the fact
that there is a Legendrian approximation $L_1$ of the binding
component $T$, which is included in a page $\Sigma$ (see
\cite{vv}). Then we express the effect of  a full Lutz twist along
$T$ by a contact $(+1)$-surgery on a four-component link
$\mathbb{L} =L_1 \sqcup L_2 \sqcup L_3 \sqcup L_4$ in $M$, where
$L_2$ is a Legendrian push-off of $L_1$ with two additional
up-zigzags, $L_3$ is a Legendrian push-off of $L_2$ and $L_4$ is a
Legendrian push-off of $L_3$ with two additional up-zigzags.  This
result is, indeed, analogues to the results in \cite{dgs2, etn}.
Next, we stabilize the open book at hand to embed all four
components of the Legendrian link $\mathbb{L}$ into a page (cf.
\cite{eobd}). Finally, we use the fact that, a contact
$(+1)$-surgery on $\mathbb{L}$ corresponds to additional
left-handed Dehn twists along each $L_i$ ($i=1,\ldots, 4)$, on the
page. As a result, we observe that $(\Sigma', \phi')$ is obtained
from $(\Sigma, \phi)$ by a \emph{local} modification near the
binding and,   by construction, the genus of $\Sigma$ is the same
as the genus of $\Sigma'$.

With an other point of view, what we describe in this note is a
\emph{planar} relative open book that is compatible with the
contact $T^2 \times I$ layer we ``insert" while performing a full
Lutz twist. Relative open books for contact $3$-manifolds whose
boundary consists of convex tori were introduced and studied by
Van Horn-Morris in his thesis \cite{morris}. In particular, he
describes a relative open book compatible with a basic slice. It
turns out that the relative open book we obtain for the full Lutz
twist through contact surgery, can be essentially obtained by
appropriately gluing together a string of relative open books for
the basic slice.

Throughout this  paper, we assume that all contact structures are
positive and co-oriented, all transverse knots are positively
transverse, and all stabilizations are positive.  The reader may
turn to \cite{et,eobd,geiges,os} for the basic material on contact
topology.

\section{Lutz twists}\label{lt}

Let $T$ be a knot transverse to the contact structure $\xi$ in a
$3$-manifold $M$.  Then, in suitable local coordinates, we can
identify $T$ with $S^1\times \{0\}\subset S^1\times D_{\delta}^2$
for some, possibly small $\delta > 0$ such that $\xi=\ker (d\theta
+ r^2 d\varphi)$ and $\partial_{\theta}$ corresponds to the
positive orientation of $T$. In order to simplify the notation, we
will work with $S^1\times D^2$ as a local model.  A simple Lutz
twist along $T$ is defined by replacing the contact structure
$\xi$ on $M$ by $ \xi^T$ which coincides with $\xi$ outside the
solid torus $S^1\times D^2$ and on $S^1\times D^2$ is given by
$$\ker (h_1(r)d\theta + h_2(r)d\varphi)$$ where $h_1, h_2
\colon [0, 1]\to \bR$ are smooth functions satisfying the following
conditions:
\begin{itemize}
\item[(i)]$h_1(r)=-1$ and $h_2(r)=-r^2$ near $r=0$,
\item[(ii)]$h_1(r)=1$ and $h_2(r)=r^2$ near $r=1$,
\item[(iii)]$(h_1(r), h_2(r))$ is never parallel to $(h'_1(r), h'_2(r))$ for $r\neq 0$.
\end{itemize}
Note that $\xi^T$ is well-defined up to isotopy, i.e., the isotopy
class of $\xi^T$ does not depend on the particular choice of the
functions $h_1$ and $h_2$. Moreover, it is clear that a simple Lutz
twist does not change the topology of the underlying $3$-manifold,
but,  in general, $\xi$ and $\xi^T$ are not homotopic as oriented
$2$-plane fields (see \cite[Section 4.3]{geiges}).

A full Lutz twist along $T$ is defined similar to a simple Lutz
twist but the boundary conditions (i) and (ii) above are replaced
by
$$ h_1(r)=1 \ \textrm{and} \ h_2(r)=r^2 \ \textrm{for} \ r\in [0,
\varepsilon]\cup[1-\varepsilon, 1]$$ for some small $\varepsilon$,
(iii) still holds and we require that the curve $(h_1(r), h_2(r))$
completes a full turn around the origin.  A full Lutz twist does
not change the homotopy class of the contact structure as a
$2$-plane field, nor the topology of the underlying manifold (see
\cite[Proposition 4.5.4]{geiges}). Let $\zeta$ denote the contact
structure obtained by applying a full Lutz twist along $T$. Note
that $(M, \xi^T)$ and $(M, \zeta)$ are both overtwisted.

\begin{remark} \label{otd}
Full Lutz twist is equivalent to replacing the contact structure
in the
 $T^2\times I$ layer (corresponding to $r \in I=[\epsilon,
1-\epsilon]$) by another contact structure with an additional
$2\pi$-twist in the $I$-direction.
\end{remark}

\section{The surgery diagram for a full Lutz twist}

In a recent series of papers \cite{dg1,dg2,dgs1}, a notion of
contact $r$-surgery along Legendrian knots in contact
$3$-manifolds is described, where $r \in (\mathbb{Q} \setminus
\{0\}) \cup \{\infty\}$ denotes the framing relative to the
natural contact framing. This generalizes the contact surgery
introduced by Eliashberg \cite{elias} and Weinstein \cite{we},
which corresponds to the contact $(-1)$-surgery.

On the other hand,  the classical notion of a Lutz twist (see
\cite{lut,mar}) played an important role in constructing various
contact structures. It turns out that, a \emph{simple} Lutz twist
along a transverse knot in a contact $3$-manifold is equivalent to
contact $(+1)$-surgery along a Legendrian two-component link
\cite{dg2}. Moreover, an explicit Legendrian surgery description
for the simple Lutz twist is given in \cite{dgs2,etn}. Similarly,
a \emph{full} Lutz twist along a transverse knot in a contact
$3$-manifold is equivalent to contact $(+1)$-surgery along a
Legendrian four-component link (cf. \cite{dg2,eh}). The idea here
is that a full Lutz twist is equivalent to two iterated simple
Lutz twists.

\begin{figure}[ht]
  \relabelbox \small {
  \centerline{\epsfbox{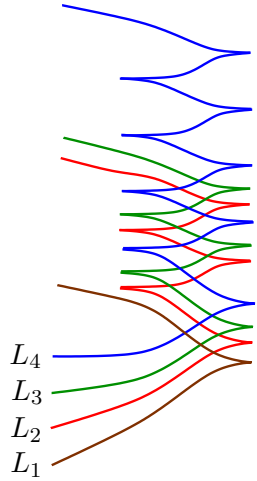}}}

\relabel{1}{{$L_1$}} \relabel{2}{{$L_2$}} \relabel{3}{{$L_3$}}
\relabel{4}{{$L_4$}}
  \endrelabelbox
        \caption{Legendrian link $\mathbb{L} =
L_1\sqcup L_2\sqcup L_3\sqcup L_4$ }
        \label{lutz}
\end{figure}

\begin{theorem} \label{alutz}
Let $L_1$ be an oriented Legendrian knot in $(M, \xi)$, represented
by its front projection in $(\bR^3, \xi_{st})$ disjoint from the
link describing $(M, \xi)$ and $L_2$ be a Legendrian push-off of
$L_1$ with two additional up-zigzags, $L_3$ be a Legendrian push-off
of $L_2$ and $L_4$ be a Legendrian push-off of $L_3$ with two
additional up-zigzags. Let $\mathbb{L}:= L_1\sqcup L_2\sqcup
L_3\sqcup L_4$ (see Figure~\ref{lutz}) and $\xi'$ be the contact
structure obtained from $\xi$ by  contact $(+1)$-surgery on
$\mathbb{L}$. If $\zeta$ denotes the contact structure obtained from
$\xi$ by a full Lutz twist along a positive transverse push-off $T$
of $L_1$, then $\xi'$ and $\zeta$ are isotopic.
\end{theorem}

\begin{proof}
We first show that contact $(+1)$-surgery on the Legendrian link
$\mathbb{L}$ does not topologically change the  underlying
manifold $M$. To see this, note that an additional zigzag adds a
negative twist to the contact framing. Hence, topologically
contact $(+1)$-surgery on $L_4$ is the same as a contact
$(-1)$-surgery along a Legendrian push-off of $L_3$. Therefore, by
\cite[Proposition 8]{dg1}, the contact $(+1)$-surgery on $L_3$
topologically  cancels out the contact $(+1)$-surgery on $L_4$.
The same argument holds for the contact $(+1)$-surgeries on $L_1$
and $L_2$.

We know that $\zeta$ is overtwisted. It is not too hard to see
that $\xi'$ is also overtwisted (cf. \cite{oz}). Once we show that
$\xi'$ is homotopic to $\zeta$ as an oriented $2$-plane field,
then the result immediately follows from Eliashberg's
classification of overtwisted contact structures \cite{eli}. Since
a full Lutz twist does not change the homotopy class of $\xi$ as a
$2$-plane field, i.e., $\xi$ is homotopic to $\zeta$, we need to
verify that $\xi$ is homotopic to $\xi'$. Recall that for any two
$2$-plane fields $\xi$ and $\xi'$ on $M$, there is an obstruction
$d^2(\xi, \xi')\in H^2(M; \mathbb{Z})$ for $\xi$ to be homotopic
to $\xi'$ over the $2$-skeleton of $M$ and if $d^2(\xi, \xi')=0$,
after applying a homotopy which takes $\xi$ to $\xi'$ over the
$2$-skeleton, there is another obstruction $d^3(\xi, \xi')$ for
$\xi$ to be homotopic to $\xi'$ over all of $M$.

Consider the standard tight contact $(S^1\times S^2, \xi)$, which
can be represented by contact $(+1)$-surgery on a Legendrian
unknot $L_0$ with only two cusps.  Let $L_1$ be a Legendrian
push-off of $L_0$. Note that, by the neighborhood theorem for
Legendrian knots, it suffices to prove the vanishing of the
two-dimensional obstruction $d^2(\xi, \xi')$ for this particular
$L_1$ (cf. \cite{dgs2}).  It is well-known that $e(\xi)=0$. Here
we claim that $e(\xi') = 0$, as well. It follows that $d^2(\xi,
\xi')=0$, by the formula $2d^2(\xi, \xi')=e(\xi)-e(\xi')$ (see
\cite[Remark 4.3.4]{geiges}).

The Thurston-Bennequin invariants of the Legendrian knots $L_0, L_1,
\ldots, L_4 $ can easily be computed from their front projections as
$tb(L_0)=-1$, $tb(L_1)=-1$, $tb(L_2)=tb(L_3)=-3$ and $tb(L_4)=-5$.
Thus, the topological framings of the surgeries are given by
$tf(L_0)=tf(L_1)=0$, $tf(L_2)=tf(L_3)=-2$ and $tf(L_4)=-4$. Write
$\mu_i$ for the meridional circle to $L_i$ as well as the homology
classes they represent in the homology of the surgered manifold. It
is well-known that $H_1(M; \mathbb{Z})$ is generated by the
meridians $\{ \mu_0, \ldots, \mu_4\}$ with relations given by
$$tf(L_i) \mu_i + \sum_{j\neq i} lk(L_i, L_j)\mu_j =0, \ i=0,
\ldots, 4 .$$ These equations imply that
$\mu_0=\mu_1=\mu_4=-\mu_2=-\mu_3$. Now with PD denoting the
Poincar\'e duality isomorphism, we have (see \cite[Corollary
3.6]{dgs1})
\begin{align*}
e(\xi')&= \Sigma_{i=1}^{4} \tr{rot}(L_i)\tr{PD}^{-1}(\mu_i)\\
&=-2\tr{PD}^{-1}(\mu_2)-2\tr{PD}^{-1}(\mu_3)-4\tr{PD}^{-1}(\mu_4)=0.
\end{align*}

Finally, let us see the effect of the surgery along $\mathbb{L}$
on the $3$-dimensional obstruction.  It is sufficient to consider
an oriented knot $L_1$ in $(S^3, \xi_{st})$.  The absolute
$d_3$-invariant (for $2$-plane fields in $S^3$) of the contact
structure $\xi'$ obtained by these surgeries is given by (see
\cite[Corollary 3.6]{dgs1})

$$d_3(\xi')
= \frac{1}{4}(c^2-3\sigma(X)-2\chi(X))+q \; ,$$

\noindent where $X$ denotes the handlebody obtained from $D^4$ by
attaching four $2$-handles corresponding to the surgeries, $q$
denotes the number of components in $\mathbb{L}$ on which we
perform $(+1)$ surgery and $c\in H^2(X; \mathbb{Z})$ is given by
$c([\Sigma_i])=rot(L_i)$ on $[\Sigma_i]\in H_2(X; \mathbb{Z})$
where $\Sigma_i$ is the Seifert surface for $L_i$. It is clear
that $\chi(X)=5$.

\begin{lemma} We have $\sigma(X)=0$ and $ c^2 = -8$.
\end{lemma}

\begin{proof} Let $t$ denote the Thurston-Bennequin invariant of $L_1$.  Hence we
have $tb(L_2)=tb(L_3)=t-2$ and $tb(L_4)=t-4$. Then the topological
framings of the surgeries are
$$tf(L_1)=t+1,\; tf(L_2)=tf(L_3)=t-1,\;  \tr{and} \  \ tf(L_4)=t-3 \
.$$  The linking number between $L_1$ and  $L_j$ is given by
$lk(L_1, L_j)=t$ for $j=2, 3$ and $4$.  Also we have $lk(L_2,
L_3)=lk(L_2, L_4)=t-2$ and $lk(L_3, L_4)=t-2$.  Then $\sigma(X)$ is
the signature of the linking  matrix
$$\left[
\begin{array}{cccc} t+1&t&t&t\\ t&t-1&t-2&t-2
\\t&t-2&t-1&t-2\\t&t-2&t-2&t-3
\end{array} \right]$$
If we slide $L_4$ over $L_3$ and slide $L_2$ and $L_3$ over $L_1$,
then the linking matrix becomes
$$A=\left[
\begin{array}{cccc} t+1&-1&-1&0\\ -1&0&-1&0
\\-1&-1&0&-1\\0&0&-1&0
\end{array} \right]$$
The characteristic polynomial for the matrix $A$ is
$\lambda^4-(t+1)\lambda^3-4\lambda^2+2(t+2)\lambda+1 .$  By
analyzing the coefficients of this polynomial  one can see that the
eigenvalues $\lambda_1, \ldots, \lambda_4$ satisfy the following
equalities:
\begin{itemize}
\item[(i)]$\lambda_1\lambda_2 +
(\lambda_1+\lambda_2)(\lambda_3+\lambda_4)+
\lambda_3\lambda_4=-4$,
\item[(ii)]$\lambda_1\lambda_2\lambda_3\lambda_4=1.$
\end{itemize}

$A$ is a real symmetric matrix, so the eigenvalues must be real
and by (ii) we have three cases for the eigenvalues of $A$:
\begin{itemize}
\item[(I)] all the eigenvalues are positive, \item[(II)] all the
eigenvalues are negative, \item[(III)] there are two positive and
two negative eigenvalues.
\end{itemize}

Clearly, (i) implies that the matrix $A$ has two positive and two
negative eigenvalues and hence $\sigma(X)=0$.

In order to compute $c^2$, set $r=rot(L_1)$. Then
$rot(L_2)=rot(L_3)=r-2,\;  \ \tr{and} \ rot(L_4)=r-4$. As in Section
$3$ of \cite{dgs1}, we have $$c^2 = xr+y(r-2)+z(r-2)+w(r-4) ,$$
where $(x, y, z, w)$ is the solution of the system of equations
$$ \left[
\begin{array}{cccc} t+1&t&t&t\\ t&t-1&t-2&t-2
\\t&t-2&t-1&t-2\\t&t-2&t-2&t-3
\end{array} \right] \left[
\begin{array}{cccc} x\\ y
\\z\\w
\end{array} \right] = \left[
\begin{array}{cccc} r\\ r-2
\\r-2\\r-4
\end{array} \right].$$ It  follows that $x=r, y=z=-2-r,
w=4+r$, and hence $c^2 = -8$.

\end{proof}

\noindent Consequently,

$$
d_3(\xi')= \frac{1}{4}(c^2-3\sigma(X)-2\chi(X))+q = -1/2 =
d_3(\xi_{st}), $$ which implies that $d^3(\xi, \xi') =0$.
Therefore, since $d^2(\xi, \xi')=d^3(\xi, \xi')=0$, we conclude
that $\xi$ is homotopic to $\xi'$, i.e.,  the contact
$(+1)$-surgery on $\mathbb{L}$ does not change the homotopy class
of the contact structure.
\end{proof}

\section{The effect of a full Lutz twist along the binding of an open book}

In this section, we describe the effect of a full Lutz twist along
the binding of an open book. Let $T$ denote a binding component of
an open book $(\Sigma, \phi)$ compatible with a closed contact
$3$-manifold $(M, \xi)$. The first step of our construction is to
realize a push-off of $T$ as a Legendrian curve on the page
$\Sigma$.  If $T$ is \emph{not} the only binding component of
$(\Sigma, \phi)$, then a push-off of $T$  can be Legendrian
realized on $\Sigma$, otherwise this is impossible. The next
result, which is due to Vela-Vick \cite{vv}, allows us to handle
the case where the binding of $(\Sigma, \phi)$ is connected.

\begin{lemma} Suppose that $T$ is the
only binding component of $(\Sigma, \phi)$. Then there is a
Legendrian approximation $L_1$ of $T$, realized as a curve on a
page of the open book $(\Sigma_1, \phi_1)$ obtained from $(\Sigma,
\phi)$ by stabilizing once, without changing the transverse
isotopy type of $T$.
\end{lemma}

\begin{proof} Although this was shown in \cite[Lemma 3.1]{vv},  here we
give an alternative proof of this fact. Suppose that $T$ is the
only binding component of $(\Sigma, \phi)$. We take a small
meridional circle $m$ linking $T$ once and hence intersecting
every page $\Sigma$ transversely. Now we perform $+1$ surgery
along $m$ and extend the open book into the glued in solid torus
as follows (see \cite{eo}): Foliate the solid torus we glue in by
annuli so that each annulus has the core of this solid torus as
one of its boundary components. The other boundary component of
each annulus is a $(1,1)$-curve which is identified with a
meridional circle to $m$ during the surgery. So we delete a disk
from every page $\Sigma$ and glue in an annulus to get an open
book with page $\Sigma_1$ and monodromy $\phi_1$, where $\phi_1$
is obtained from $\phi$ by composing with a right-handed Dehn
twist along a curve parallel to the newly introduced binding
component. In fact, we just stabilized $(\Sigma, \phi)$. The
contact structure compatible with $(\Sigma_1, \phi_1)$ is isotopic
to $\xi$ and we did not touch the transverse knot $T$ in $M$ while
performing the surgery.

\begin{figure}[ht]
  \relabelbox \small {
  \centerline{\epsfbox{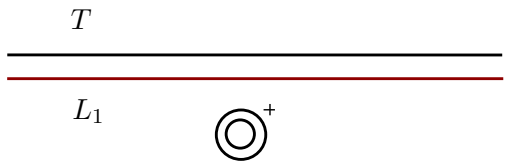}}}

\relabel{1}{{$T$}} \relabel{2}{{$L_1$}}
  \endrelabelbox
        \caption{Legendrian knot $L_1$ on the page $\Sigma_1$}
        \label{lutzbinding}
\end{figure}

The knot  $L_1$ which is a push-off of $T$ on the new page
$\Sigma_1$ as shown abstractly in Figure~\ref{lutzbinding} can now
be Legendrian realized since it is homologically nontrivial on
$\Sigma_1$. Consequently, $L_1$ is a Legendrian knot on the page
$\Sigma_1$ whose positive transverse push-off is $T$.

\end{proof}

\begin{figure}[ht]
  \relabelbox \small {
  \centerline{\epsfbox{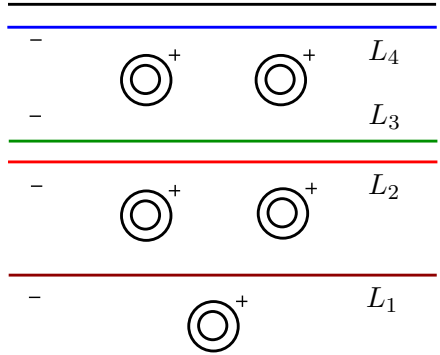}}}

\relabel{1}{{$L_1$}} \relabel{2}{{$L_2$}} \relabel{3}{{$L_3$}}
\relabel{4}{{$L_4$}}
  \endrelabelbox
        \caption{Modification near the binding which corresponds to the effect of a full Lutz twist.
        The solid curves without a sign are the binding components.}
        \label{lutzop}
\end{figure}

Since $L_2$ is obtained from a push-off of $L_1$ by adding two
zigzags, we can realize $L_2$ on a page of an open book
$(\Sigma_2, \phi_2)$ obtained by stabilizing $(\Sigma_1, \phi_1)$
twice (cf. \cite{etn}).  To be more precise, $L_2$ is a push-off
of $L_1$ on $\Sigma_2$, except that $L_2$ goes over the two new
$1$-handles glued to $\Sigma_1$ in the stabilization process. By
continuing in this manner, we see that there is an open book
$(\Sigma, \phi)$, compatible with $(M, \xi)$, containing the
Legendrian link $\mathbb{L}$ on a page. Then the open book
$(\Sigma, \phi \circ D_\mathbb{L}^-)$ is compatible with $(M,
\xi')$, where $D_\mathbb{L}^-$ denote the composition of
left-handed Dehn twists along each component of the link
$\mathbb{L} \subset \Sigma$ (see Figure~\ref{lutzop}).
Consequently, by Giroux's correspondence \cite{giroux} coupled
with Theorem~\ref{alutz},  we conclude that $(\Sigma, \phi \circ
D_\mathbb{L}^-)$ is compatible with $(M, \zeta)$.

\begin{figure}[ht]
  \relabelbox \small {
  \centerline{\epsfbox{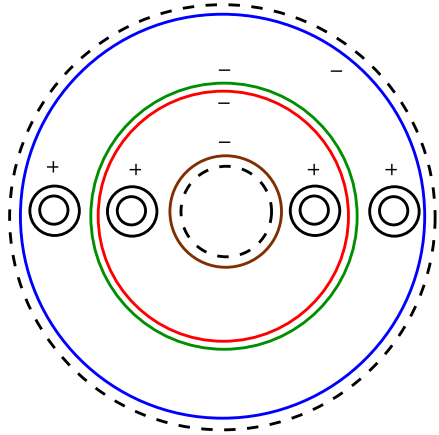}}}

  \endrelabelbox
\caption{Relative open book compatible with the contact $T^2
\times I$ layer with a $2\pi$-twist in the $I$-direction.  The
dashed curves induce the boundary $\partial (T^2 \times I)$. }
        \label{lu}
\end{figure}

\begin{remark} \label{otd}
The discussion above gives an explicit \emph{relative} open book
 for the full Lutz twist. In other words, the
contact $T^2 \times I$ layer with a $2\pi$-twist inserted while
performing a full Lutz twist has a compatible relative open book
whose page is a six-punctured sphere, and whose monodromy is the
product of the Dehn twists along the curves shown in
Figure~\ref{lu}. Note that a relative open book compatible with a
basic slice is described in Section 4.2 in \cite{morris}. It turns
out that the relative open book we obtain for the full Lutz twist
through contact surgery, can be essentially obtained by
appropriately gluing together a string of relative open books
compatible with the basic slice.
\end{remark}

\vspace{0.2in} \noindent{\bf {Acknowledgements.}} We are grateful
to Hansj\"{o}rg Geiges and John B. Etnyre for helpful
conversations. The authors were partially supported by the
research grant 107T053 of the Scientific and Technological
Research Council of Turkey.

\bibliographystyle{amsplain}
\providecommand{\bysame}{\leavevmode\hbox
to3em{\hrulefill}\thinspace}
\providecommand{\MR}{\relax\ifhmode\unskip\space\fi MR }

\providecommand{\MRhref}[2]{%
  \href{http://www.ams.org/mathscinet-getitem?mr=#1}{#2}
} \providecommand{\href}[2]{#2}

\end{document}